\documentclass[12pt]{amsart}
\usepackage{amsmath,amssymb,amsfonts,amsthm,amsopn}
\usepackage{graphics}

\newcommand{\ft}{Fourier transform}
\newcommand{\stft}{short-time Fourier transform}

\newcommand{\tf}{time-frequency}

\newcommand{\modsp}{modulation space}

\newtheorem{tm}{Theorem}[section]
\newtheorem{lemma}[tm]{Lemma}

\newtheorem{theorem}{Theorem}[section]

\newtheorem{example}[theorem]{Example}

\newtheorem{proposition}[theorem]{Proposition}
\newtheorem{remark}[theorem]{Remark}

\newcommand{\beqa}{\begin{eqnarray*}}
\newcommand{\eeqa}{\end{eqnarray*}}

\newcommand{\field}[1]{\mathbb{#1}}
\newcommand{\bR}{\field{R}}        
\newcommand{\bN}{\field{N}}        
\newcommand{\bZ}{\field{Z}}        
\def\la{\lambda}

 \def\cF{\mathcal{F}}              
 \def\cS{\mathcal{S}}
 \def\cD{\mathcal{D}}

 \def\cG{\mathcal{G}}
 \def\cM{\mathcal{M}}

 \def\cC{\mathcal{C}}

 \def\cX{\mathcal{X}}

\def\a{\aleph}

\def\vgf{V_gf}

\def\rd{\bR^d}

\def\rdd{{\bR^{2d}}}
\def\zdd{{\bZ^{2d}}}

\def\lrd{L^2(\rd)}
\def\lrdd{L^2(\rdd)}
\def\zd{\bZ^d}

\def\intrd{\int_{\rd}}
\def\intrdd{\int_{\rdd}}

\def\R{\right)}

\def\<{\left<}
\def\>{\right>}

\def\mv1{M_v^1}

\def\mpq{M^{p,q}}
\def\Mmpq{M_\mu^{p,q}}
\def\phas{(x,\o )}
\def\mn{(m,n)}
\def\mn'{(m',n')}

\def\o{\eta}
\def\a{\alpha}
\def\b{\beta}
\def\z{\zeta}

\def\ZZ{\mathbb{Z}}

\def\N{\mathbb{N}}
\def\R{\mathbb{R}}
\def\Ren{\mathbb{R}^d}
\def\Renn{\mathbb{R}^{2d}}

\def\sch{\mathcal{S}}

\def\Fur{\mathcal{F}}

\def\Sn2{S_{2}(L^{2}(\Ren))}
\def\S1{S_{1}(L^{2}(\Ren))}
\def\sig00{\sigma_{0,0}}

\def\la{\langle}
\def\ra{\rangle}

\begin{document}

\begin{abstract} We use \tf\, methods for the study of Fourier Integral operators (FIOs). In this paper we shall show  that Gabor frames
provide very efficient representations for a large class of FIOs.
Indeed, similarly to the case of shearlets and curvelets frames \cite{candes,Guo-labate}, the matrix representation of a Fourier Integral Operator with respect to a Gabor frame is well-organized.  This is used as  a powerful tool to
 study the boundedness of FIOs on \modsp s.  As special cases,  we   recapture boundedness results on modulation spaces for pseudo-differential operators   with symbols in  $M^{\infty,1}$ \cite{grochenig-heil},
 for some unimodular Fourier multipliers \cite{benyi}
 and metaplectic operators \cite{cordero2,book}.
\end{abstract}

\title{  Time-Frequency Analysis of  Fourier Integral Operators}
\author{Elena Cordero, Fabio Nicola and Luigi Rodino}
\address{Department of Mathematics,
University of Torino, via
Carlo Alberto 10, 10123
Torino, Italy}
\address{Dipartimento di Matematica,
Politecnico di Torino, corso
Duca degli Abruzzi 24, 10129
Torino, Italy}
\address{Department of Mathematics,
University of Torino, via
Carlo Alberto 10, 10123
Torino, Italy}
\email{elena.cordero@unito.it}
\email{fabio.nicola@polito.it}
\email{luigi.rodino@unito.it}
\thanks{}
\maketitle

\section{Introduction}

Fourier Integral Operators (FIOs)  are a mathematical
tool to study variety of problems arising in partial
differential equations. Originally introduced by Lax
\cite{Lax} for the construction of parametrices in the
Cauchy problem for  hyperbolic equations, they have
 been widely employed to represent  solutions  to
 Cauchy problems, in the framework of  both pure and
  applied mathematics (see,
  e.g.,the papers
   \cite{boulkhemair,candes,duistermat-guillemin,duistermat-hormander,Guo-labate,hormander0},
   the books
   \cite{hormander,stein93,treves}
and
   references therein).
In particular, they were employed by Helffer and Robert
 \cite{helffer84,helffer-rob1} to study the spectral property of a class of globally elliptic operators, generalizing the harmonic oscillator of
 the Quantum Mechanics. The Fourier Integral operators we work with, possess a phase function similar to  those of \cite{helffer84,helffer-rob1}. A simple example is the resolvent of   the Cauchy problem
  for the Schr\"odinger
  equation with a quadratic
  Hamiltonian.

  For a given function $f$ on $\rd$  the
\emph{Fourier Integral Operator} (FIO) $T$ with symbol $\sigma$  and phase
$\Phi$ on $\rdd$ can be formally defined by
\begin{equation*}
    Tf(x)=\intrd e^{2\pi i \Phi\phas} \sigma\phas \hat{f}(\o)d\o.
\end{equation*}
The phase
function $\Phi(x,\eta)$ is smooth on $\rdd$, fulfills the  estimates
\begin{equation}\label{A}
|\partial_z^\a \Phi(z)|\leq
C_\a,\quad |\a|\geq 2,\quad
z\in\rdd,
\end{equation}
and the nondegeneracy condition
\begin{equation}\label{B}
    |\det\,\partial^2_{x,\eta} \Phi(x,\o)|\geq \delta>0,\quad \quad (x,\o)\in\rdd.
\end{equation}
The symbol $\sigma$ on $\rdd$
satisfies
\begin{equation*}
    |\partial_z^\a \sigma(z)|\leq
    C_{\a},\quad\ |\a|\leq N,\quad
    \mbox{a.e.}\, z\in\rdd,
\end{equation*}
for a fixed $N>0$
(in the sequel we shall work also with rougher symbols).

The first goal of this paper
is to rephrase the operator
$T$ in terms of
time-frequency analysis (see
Gr\"ochenig \cite{book} and
the next Section $2$ for a
review of the time-frequency
methods.) Denoting
$T_{x}f(t)=f(t-x)$, $M_\o
f(t)=e^{2\pi i \o t} f(t)$,
for $\a,\beta>0$, $g\in\lrd$,
the set of time-frequency
shifts
$\cG(g,\a,\b)=\{g_{m,n}:=M_{
n}T_{ m} g\}$ with
$(m,n)\in\a \zd\times\b\zd $,
is a \emph{ Gabor frame} if
there exist positive
constants $A,B>0$, such that
\begin{equation}\label{freme}
A \|f\|_{L^2}\leq \sum_{m,n}
|\la f,T_{ m}M_{ n}
g\ra|^2\leq B
\|f\|_{L^2},\quad \forall f
\in\lrd.
 \end{equation} In
Section $3$ we show that the
matrix representation of a
FIO $T$ with respect to a
Gabor frame with
$g\in\mathcal{S}(\R^d)$ is
well-organized (similarly to
frames of curvelets and
shearlets
\cite{candes,Guo-labate}),
provided that the symbol
$\sigma$ satisfies the decay
estimate for every $N>0$ (see
Theorem \ref{well-pose}):
\vskip0.3truecm
\textbf{Theorem 1.} \emph{For
each $N>0$, there exists a
constant $C_N>0$ such that
\begin{equation}\label{C}
    |\la T g_{m,n}, g_{m',n'}\ra|\leq C_N {\la \chi(m,n)
    -(m',n')\ra^{-2N}
    },
\end{equation}
where  $\chi$ is the canonical
transformation generated by $\Phi$. }
\vskip0.3truecm In the special case of
pseudodifferential operators such an
almost diagonalization was already
obtained in
\cite{grochenig-ibero,rochberg}.
Indeed, notice that pseudodifferential
operators correspond to the phase
$\Phi(x,\eta)=x\eta$ and canonical
transformation
$\chi(y,\eta)=(y,\eta)$.\par As a
rilevant byproduct of the results of
Section $3$, we study the boundedness
properties of the operator $T$ on the
so-called modulation spaces (Section
$4$ and $5$). To define them, we fix a
non-zero Schwartz function $g$ and
consider the short-time Fourier
Transform $V_g f$ of a function $f$ on
$\rd$ with respect to $g$
$$ V_gf(x,\o)=\la f,M_\o T_x g\ra =\int_{\Ren}
 f(t)\, {\overline {g(t-x)}} \, e^{-2\pi i\o t}\,dt\, $$
 which provides a time-frequency representation of  $f$. The (unweighted) modulation space $M^{p,q}$ is the closure of the Schwartz class with respect to the norm
 $$\|f\|_{M^{p,q}}=\|V_gf\|_{L^{p,q}}=\left(\int_{\Ren}
  \left(\int_{\Ren}|V_gf(x,\o)|^p\,
    dx\right)^{q/p}d\o\right)^{1/p}  \,
 $$
(with appropriate modifications when $p=\infty$ or $q=\infty$). In particular, when $p=q$ we simply write $M^{p,p}=M^{p}$, see Subsection \ref{modsp} for  exhaustive definitions and properties.

These spaces  were introduced by Feichtinger \cite{F1} and have become canonical for both time-frequency and phase-space analysis \cite{F06},  most recent employment being  the study of PDEs \cite{benyi,
benyi2,baoxiang3,baoxiang2,baoxiang}.
\par
If $g\in M^1$, and
 the Gabor frame $\{T_{ m}M_{ n} g
 ;
  (m,n)\in\a \zd\times\b\zd  \}$  is a {\it tight frame},
   namely \eqref{freme} holds with $A=B$,
  then it extends to a  Banach frame for the
   \modsp s $\mpq (\rdd )$, with the norm equivalence
$$ \| f\|_{M^{p,q}}\asymp
 \|\la f ,T_{ m}M_{ n} g\ra _{m,n}\|_{l^{p,q}}.
$$
Then,  boundedness of the FIO $T$ on $M^{p,q}$ is equivalent  to that of the infinite matrix $\la T g_{m,n},\, g_{m',n'}\ra$ on the spaces of sequences $l^{p,q}$.

Whence, the estimates \eqref{C}  readily give (see  Theorem \ref{contmp} for a more general version):
\vskip0.3truecm
\textbf{Theorem 2.} \emph{For $N>d$, $1\leq p<\infty$, the Fourier integral operator $T$, with symbol $\sigma$ and phase $\Phi$ as above, extends to a continuous operator on $M^p$.\\
 \noindent(In the case $p=\infty$, the space $M^\infty$ is replaced by the closure of the Schwartz function with respect to $\|\cdot\|_{M^\infty}$).}
\vskip0.3truecm

The continuity property of a FIO  $T$ on $M^{p,q}$,
 with $p\not=q$, fails  in general. Indeed,  an example
  is  provided by the operator
  $T f(x)=e^{\pi i|x|^2}f(x)$, corresponding to
  $\Phi(x,\eta)=x\eta+|x|^2/2$, $\sigma\equiv1$,
  which is bounded
  on $M^{p,q}$ if and only if $p=q$
  (see Proposition \ref{contro}).

Hence, we introduce a new
condition on the phase
$\Phi$, namely that
  the map $x\longmapsto\nabla_x\Phi(x,\eta)$
 has a range of finite diameter, uniformly
with respect to $\eta$, that
allows us to get the
boundedness on $M^{p,q}$
(Theorem \ref{cont}):
\vskip0.3truecm
\textbf{Theorem 3.} \emph{For
$N>d$, $1\leq p, q<\infty$,
the Fourier integral operator
$T$, with symbol $\sigma$ and
phase $\Phi$ as above, and
such that
\begin{equation*}
\sup_{x,x',\eta}\left|
\nabla_x\Phi(x,\eta)-\nabla_x\Phi(x',\eta)\right|<\infty,
\end{equation*}
extends to a continuous operator on $M^{p,q}$.\\
\noindent (In the case $p=\infty$ or
$q=\infty$, the space $M^{p,q}$ is
replaced by the closure of the Schwartz
function with respect to
$\|\cdot\|_{M^{p,q}}$).}
\vskip0.3truecm As a particular case,
we recapture recent boundedness results
of unimodular Fourier multipliers
\cite{benyi} (see Example
\ref{osser3}).\\ With respect to
Theorem 2, here the proof is
 more delicate and combines the estimate \eqref{C}
 with a generalized version of Schur's Test (Proposition \ref{proschur}).

To have a simple idea of the possible applications of Theorems $1$, $2$ and $3$,   consider the Cauchy problem
  \begin{equation}\label{C1}
\begin{cases} i \displaystyle\frac{\partial
u}{\partial t} +H u=0\\
u(0,x)=u_0(x),
\end{cases}
\end{equation}
where $H$ is the Weyl quantization of a quadratic form on
$\R^d\times \R^d$ (see, e.g., \cite{cordero2,folland}). Simple examples are $H=-\frac{1}{4\pi}\Delta+\pi|x|^2$, or
$H=-\frac{1}{4\pi}\Delta-\pi|x|^2$ (see \cite{berezinshubin}). The solution to \eqref{C1}  is a one-parameter family of FIOs:  $$u(t,x)=e^{itH}u_0, $$ with symbol $\sigma\equiv 1$ and a phase given by a quadratic form $\Phi\phas$, satisfying trivially the preceding assumptions \eqref{A} and \eqref{B} (see \cite{folland} for details).  We address to Section $7$ for a debited study of such operators.

Finally, Section $6$ presents a variant
of Theorem $1$,  cf.  \eqref{symbstft},
and a generalization of Theorem $2$,
cf.  Theorem \ref{contmp2},  to the
case of FIOs $T$ with symbols in the
modulation space $M^{\infty,1}$. This
generalizes the known boundedness
results on $M^p$ of pseudodifferential
operators with symbols in
$M^{\infty,1}$ \cite{grochenig-heil},
and intersects a previous result of
Boulkhemair \cite{boulkhemair}  on
$L^2$ boundedness of FIOs. We address
also to recent contribution
\cite{concetti-toft}, where the
continuity and Schatten-von Neumann
properties of similar operators when
acting on $L^2$ are proved. \vspace{5
mm}

\textbf{Notation.} We define
$|t|^2=t\cdot t$, for $t\in\Ren$, and
$xy=x\cdot y$ is the scalar product on
$\Ren$.

The Schwartz class is denoted by
$\sch(\Ren)$, the space of tempered
distributions by  $\sch'(\Ren)$.   We
use the brackets  $\la f,g\ra$ to
denote the extension to $\sch
(\Ren)\times\sch '(\Ren)$ of the inner
product $\la f,g\ra=\int f(t){\overline
{g(t)}}dt$ on $L^2(\Ren)$. The Fourier
transform is normalized to be ${\hat
  {f}}(\o)=\Fur f(\o)=\int
f(t)e^{-2\pi i t\o}dt$, the involution
$g^*$ is $g^*(t) = \overline{g(-t) }$
and the inverse Fourier transform is
${\check  f}(\o)=\Fur^{-1}f (\o)={\hat
{f}}(-\o)$.\par Translation and
modulation ({\it time and frequency
shifts}) are defined, respectively, by
$$ T_xf(t)=f(t-x)\quad{\rm and}\quad M_{\o}f(t)= e^{2\pi i \o  t}f(t).$$
We have the formulas $(T_xf)\hat{} =
M_{-x}{\hat {f}}$,  $(M_{\o}f)\hat{}
=T_{\o}{\hat {f}}$, and
$M_{\o}T_x=e^{2\pi i x\o}T_xM_{\o}$.
For $\a=(\a_1,\dots,\a_d),
\a=(\a_1,\dots,\a_d) \in \zd_+$, recall
the multi-index notation $D^\a$ and
$X^\b$ for the operators of
differentiation and multiplication
$$D^\a f=\prod_{j=1}^d\partial^{\a_j}_{t_j} f\quad\quad\mbox{and}\quad\quad X^\b f(t)=\prod_{j=1}^d t_j^{\b_j} f(t),
$$
where $t=(t_1,\dots,t_d)$. We write  $dx\wedge d\xi=\sum_{j=1}^d dx_j\wedge
  d\xi_j$
 for the canonical symplectic
 2-form.\par
The spaces $l^{p,q}_\mu=l^ql^p_\mu$,
with weight $\mu$, are the Banach
spaces of sequences $\{a_{m,n}\}_{m,n}$
on some lattice, such that
$$\|a_{m,n}\|_{l^{p,q}_\mu}:=\left(\sum_{n}\left(\sum_{m}
|a_{m,n}|^p\mu(m,n)^p\right)^{q/p}\right)^{1/q}<\infty
$$
(with obvious changes when $p=\infty$ or $q=\infty)$.

We denote by $c_0$ the space of sequences vanishing at infinity.
Throughout the paper, we
shall use the notation
$A\lesssim B$ to indicate
$A\leq c B$ for a suitable
constant $c>0$, whereas $A
\asymp B$ if $c^{-1}B\leq
A\leq c B$ for a suitable
$c>0$.

\section{Time-Frequency Methods}

\subsection{ Short-Time Fourier Transform (STFT)}

The \stft\   (STFT)  of a distribution
$f\in\sch'(\Ren)$ with respect to a
non-zero window $g\in\sch(\Ren)$ is
$$ V_gf(x,\o)=\la f,M_\o T_x g\ra =\int_{\Ren}
 f(t)\, {\overline {g(t-x)}} \, e^{-2\pi i\o t}\,dt\, .$$

The STFT  $\vgf $ is defined
on many pairs of Banach
spaces. For instance, it maps
$L^2(\rd ) \times L^2(\rd )$
into $\lrdd $ and
$\sch(\Ren)\times\sch(\Ren)$
into $\sch(\Renn)$.
Furthermore, it  can be
extended  to a map from
$\sch'(\Ren)\times\sch'(\Ren)$
into $\sch'(\Renn)$.

We now recall the following
inequality \cite[Lemma
11.3.3]{book}, which is
useful when one needs to
change windows.
\begin{lemma}\label{changewind}
Let $g_0,g_1,\gamma\in\cS(\rd)$ such
that $\la \gamma, g_1\ra\not=0$ and let
$f\in\cS'(\rd)$. Then,
$$|V_{g_0} f\phas|\leq\frac1{|\la\gamma,g_1\ra|}(|V_{g_1} f|\ast|V_{g_0}\gamma|)\phas,
$$
for all $\phas\in\rdd$.
\end{lemma}

\subsection{Modulation Spaces }\label{modsp}
\label{modspdef} The modulation space
norms are a  measure of the joint
time-frequency distribution of $f\in
\sch '$. For their basic properties we
refer, for instance, to
\cite[Ch.~11-13]{book} and the original
literature quoted there.

For the quantitative description of
decay  properties, we use
 weight
functions  on the \tf\ plane.
In the sequel $v$ will always
be a continuous, positive,
even, submultiplicative
weight function (in short, a
submultiplicative weight),
hence $v(0)=1$, $v(z) =
v(-z)$, and $ v(z_1+z_2)\leq
v(z_1)v(z_2)$, for all $z,
z_1,z_2\in\Renn.$ A positive,
 weight function $\mu$ on
$\Renn$ belongs to $\cM_v$,
that is, is   {\it
  v-moderate} if
$ \mu(z_1+z_2)\leq Cv(z_1)\mu(z_2)$
for all $z_1,z_2\in\Renn.$

 For our investigation of FIOs we will mostly  use the polynomial weights   defined by
\begin{equation*}
  v_s(z)= v_s \phas = \la z\ra^s=(1+|x|^2+|\o|^2)^{s/2},\quad
   z=(x,\o)\in\Renn. \, \label{eqc1}\\
\end{equation*}
Given a non-zero window
$g\in\sch(\Ren)$, $\mu\in\cM_v$, and
$1\leq p,q\leq \infty$, the {\it
  modulation space} $M^{p,q}_\mu(\Ren)$ consists of all tempered
distributions $f\in\sch'(\Ren)$ such
that $V_gf\in L^{p,q}_\mu(\Renn )$
(weighted mixed-norm spaces). The norm
on $M^{p,q}_\mu$ is
$$
\|f\|_{M^{p,q}_\mu}=\|V_gf\|_{L^{p,q}_\mu}=\left(\int_{\Ren}
  \left(\int_{\Ren}|V_gf(x,\o)|^p\mu(x,\o)^p\,
    dx\right)^{q/p}d\o\right)^{1/p}  \,
$$
(with obvious changes when $p=\infty$ or $q=\infty$).
If $p=q$, we write $M^p_\mu$
instead of $M^{p,p}_\mu$, and
if $\mu(z)\equiv 1$ on
$\Renn$, then we write
$M^{p,q}$ and $M^p$ for
$M^{p,q}_\mu$ and
$M^{p,p}_\mu$ respectively.

Then  $\Mmpq (\Ren )$ is a
Banach space whose definition
is independent of the choice
of the window $g$. Moreover,
if $\mu\in\cM_v$ and $g \in
M^1_{v} \setminus \{0\}$,
then $\|V_gf
\|_{L^{p,q}_\mu}$ is an
equivalent norm for
$M^{p,q}_\mu(\Ren)$ (see
\cite[Thm.~11.3.7]{book}):
$$
\|f\|_{\Mmpq } \asymp
\|V_{g}f \|_{L^{p,q}_\mu }.
$$
\subsection{Wiener amalgam spaces}
For a detailed treatment we
refer to
\cite{feichtinger83,feichtinger80,feichtinger90,feichtinger-zimmermann98,
    fournier-stewart85}.\par

Let $g \in \cD (\Renn )$ be a test
function that satisfies $\sum_{(k,l)\in
\ZZ^{2d}}T_{(k,l)}g\equiv 1$. Let
$X(\Renn)$ be a
 Banach space of functions invariant under translations  and  with the property
 that $\cD \cdot  X\subset X$, e.g., $L^p, \cF L^p$, or $ L^{p,q}$. Then
the {\it Wiener amalgam space}
$W(X,L^{p,q}_\mu  )$ with local
component $X$ and global component
$L^{p,q}_\mu$ is defined as the space
of all functions or distributions for
which the norm
$$
\|f\|_{W(X,L^{p,q}_\mu)}=\Big( \int
_{\Ren } \big( \int _{\Ren } \|f\cdot
T_{(z_1,z_2)}g\|_{X})^p \,
\mu(z_1,z_2)^p \, dz_1 \big)^{q/p} \,
dz_2 \Big)^{1/q}
$$
is finite. Equivalently, $f \in
W(X,L^{p,q}_\mu )$ if and only if
$$\Big( \sum _{l\in \zd } \big(\sum _{k\in \zd }
\|f\cdot T_{(k,l)}g\|_{X} ^p \mu(k,l)^p
\big)^{q/p} \Big)^{1/q} < \infty \, .$$
It can be shown that  different choices
of $g\in \cD$  generate the same space
and yield equivalent norms. In the
sequel we shall use the
\emph{inclusions relations} between
Wiener amalgam spaces: if $B_1
\hookrightarrow B_2$ and $C_1
\hookrightarrow C_2$,
   \begin{equation*}
   W(B_1,C_1)\hookrightarrow W(B_2,C_2).
  \end{equation*}

We now recall the following
regularity property of the
STFT \cite[Lemma 4.1]{CG02}:

\begin{lemma}\label{amalg}
Let  $1\leq p,q\leq\infty$,
$\mu\in\cM_v$. If  $f\in
M^{p,q}_\mu(\Ren)$ and   $g\in
M^{1}_v(\Ren)$, then  $V_gf \in  W(\Fur
L^1,L^{p,q}_\mu)(\Renn)$ with norm
estimate
\begin{equation}\label{amalgeq}
\| V_gf\|_{W(\Fur
L^1,L^{p,q}_\mu)}\lesssim \|f\|_{
  M^{p,q}_\mu}\|g\|_{M^{1}_v}\, .
\end{equation}
\end{lemma}

We also give a slight
generalization of
\cite[Proposition
11.1.4]{book} and its
subsequent Remark.

\begin{proposition}\label{samp} Let $\cX$ be a separated
 sampling set in $\rdd$, that is,  there exists $\delta>0$,
  such that $\inf_{x,y\in\cX: x\not=y}|x-y|\geq \delta$.
  Then there exists a constant $C>0$ such that, if
  $F\in W(L^\infty, L^{p,q}_\mu)$ is any function
    everywhere defined on $\rdd$ and lower semi-continuous,
    then the restriction $F_{|_{\cX}}$ is
     in $\ell^{p,q}_{\tilde{\mu}}$, where
     $\tilde{\mu}=\mu_{|_{\cX}}$,
and
$$\|F_{|_{\cX}}\|_{\ell^{p,q}_{\tilde{\mu}}}\leq C \|F\|_{W(L^\infty,L^{p,q}_\mu)}.
$$
\end{proposition}
\begin{proof} One uses the arguments of
 \cite[Proposition 11.1.4]{book} and its subsequent Remark.
  We just shall highlight the
key points that make those
arguments to work under our
assumptions.\par First of
all, if $(r,s)\in\zdd$ and
$\cX$ is separated, then the
number of sampling points of
$\cX$ in $(r,s)+[0,1]^{2d}$
is bounded independently of
$(r,s)$.\par Secondly, for
$x\in\cX$ such that $x\in
(r,s)+[0,1]^{2d}$,
$$|F(x)|\mu(x)\leq C\|F \cdot T_{(r,s)}\chi_{[0,1]^{2d}}
\|_{L^\infty} \mu(r,s),
$$
Indeed, since $F$ is
everywhere defined on $\rdd$
and lower semi-continuous we
have $\sup |f|={\rm
ess\,sup}|f|$ on every box,
whereas $\mu(x)\leq C
\mu(r,s)$ because $\mu$ is
$v$-moderate and $v$ is
bounded on $[0,1]^{2d}$ (cf.
\cite[Lemma 11.1.1]{book}).
\end{proof}
\subsection{Gabor frames}
Fix a function $g\in\lrd$ and a lattice
$\Lambda =\a \zd\times \b\zd$, for
$\a,\b>0$. For $(m,n)\in \Lambda$,
define $g_{m,n}:=M_{ n}T_{m}g$. The set
of time-frequency shifts
$\cG(g,\a,\b)=\{g_{m,n}, (m,n)\in
\Lambda\}$ is called Gabor system.
Associated to $\cG(g,\a,\b)$ we define
the coefficient operator $C_g$, which
maps functions to sequences as follows:
\begin{equation}\label{analop}
    (C_gf)_{m,n}=(C_g^{\a,\b}f)_{m,n}:=\la
    f,g_{m,n}\ra,\quad (m,n)\in \Lambda,
\end{equation}
the synthesis operator
\[
D_g{c}=D_g^{\a,\b} c
=\sum_{(m,n)\in \Lambda}
c_{m,n} T_{m} M_{n}g,\quad
c=\{c_{m,n}\}_{(m,n)\in
\Lambda}
\]
and the Gabor frame operator
\begin{equation}\label{Gaborop}
    S_g f=S_g^{\a,\b}f:=D_g S_g f=\sum_{(m,n)\in \Lambda}\la f,g_{m,n}\ra g_{m,n}.
\end{equation}

The set $\cG(g,\a,\b)$  is called a
Gabor frame for the Hilbert space
$\lrd$ if $S_g$ is a bounded and
 invertible operator on $\lrd$. Equivalently, $C_g$ is bounded from $\lrd$ to $l^2(\a\zd\times\b\zd)$ with closed
  range, i.e., $\|f\|_{L^2}\asymp\|C_g f\|_{l^2}.$ If $\cG(g,\a,\b)$ is a Gabor frame for $\lrd$,
  then the so-called \emph{dual window}
   $\gamma=S_g^{-1} g$ is well-defined and the
   set $\cG(\gamma,\a,\b)$  is a frame (the so-called
   canonical dual frame of $\cG(g,\a,\b)$). Every $f\in \lrd$ posseses the frame expansion
\begin{equation}\label{frame}
    f= \sum_{(m,n)\in \Lambda}\la f,g_{m,n}\ra \gamma_{m,n}= \sum_{(m,n)\in \Lambda}\la f,\gamma_{m,n}\ra g_{m,n}
\end{equation}
with unconditional convergence in
$\lrd$, and norm equivalence:
$$\|f\|_{L^2}\asymp \|C_g f\|_{l^2}\asymp \|C_\gamma f\|_{l^2}.
$$
This result is contained in
\cite[Proposition 5.2.1]{book}. In
particular, if $\gamma=g$ and
$\|g\|_{L^2}=1$ the frame is called
\emph{normalized tight} Gabor frame and
the expansion \eqref{frame} reduces to
\begin{equation}\label{parsevalframe}
    f= \sum_{(m,n)\in \Lambda}\la f,g_{m,n}\ra g_{m,n}.
\end{equation}
If we ask for more regularity on the
window $g$, then the previous result
can be extended to suitable Banach
spaces, as shown below
\cite{fg97jfa,GL01}.
\begin{theorem}\label{teomod}
Let $\mu\in\cM_v$,
$\cG(g,\a,\b)$ be a
normalized tight Gabor frame
for $\lrd$,  with lattice
$\Lambda=\a\zd\times\b \zd$,
 and  $g\in M^1_v$. Define
 $\tilde{\mu}=\mu_{|_{\Lambda}}$.\\
(i) For every $1\leq
p,q\leq\infty$, $C_g:
M^{p,q}_\mu\to
l^{p,q}_{\tilde{\mu}}$ and
$D_g:
l^{p,q}_{\tilde{\mu}}\to
M^{p,q}_\mu$ countinuously
and, if $f\in\mpq_\mu,$ then
the Gabor expansions
\eqref{parsevalframe}
converge unconditionally in
$\mpq_\mu$ for $1\leq
p,q<\infty$
and all weight $\mu$, and weak$^\ast$-$M^\infty_\mu$ unconditionally if $p=\infty$ or $q=\infty$.\\
(ii) The following norms are equivalent on $\mpq_\mu$:\\
\begin{equation}\label{framexp}
\|f\|_{\mpq_\mu}\asymp \|C_g
f \|_{l^{p,q}_{\tilde
{\mu}}}.
\end{equation}
\end{theorem}
We also establish the
following properties. Denote
by $\tilde{{M}}^{p,q}_\mu$
the closure of the Schwartz
class in $M^{p,q}_\mu$.
Hence,
$\tilde{{M}}^{p,q}_\mu=M^{p,q}_\mu$
if $p<\infty$ and $q<\infty$.
Also, denote by
$\tilde{l}^{p,q}_{\tilde{\mu}}$
the closure of the space of
eventually zero sequences in
$l^{p,q}_{\tilde{\mu}}$.
Hence
$\tilde{l}^{p,q}_{\tilde{\mu}}={l}^{p,q}_{\tilde{\mu}}$
if $p<\infty$ and $q<\infty$.
\begin{theorem}\label{teomod2}
Under the assumptions of Theorem
\ref{teomod}, for every $1\leq
p,q\leq\infty$ the operator $C_g$ is
continuous from $\tilde{{M}}^{p,q}_\mu$
into $\tilde{l}^{p,q}_{\tilde{\mu}}$,
whereas the operator $D_g$ is
continuous from
$\tilde{l}^{p,q}_{\tilde{\mu}}$ into
$\tilde{{M}}^{p,q}_\mu$.
\end{theorem}
\begin{proof}
Since $C_g$ is continuous
from ${{M}}^{p,q}_\mu$ into
${l}^{p,q}_{\tilde{\mu}}$ it
suffices to verify that, if
$f$ is a Schwartz function
then $C_g(f)\in
\tilde{l}^{p,q}_{\tilde{\mu}}$.
This follows from the fact
that $C_g(f)\in
l^1_{\tilde{\mu}}$.
Similarly, for $D_g$ it
suffices to verify that, if
$c$ is any eventually zero
sequence, then $D_g(c)\in
\tilde{{M}}^{p,q}_\mu$. This
is true because $D_g(c)\in
M^1_\mu$.
\end{proof}
\section{Almost diagonalization of FIOs}
For a given function $f$ on $\rd$  the
FIO $T$ with symbol $\sigma$ and phase
$\Phi$ can be formally defined by
\begin{equation}\label{fio}
    Tf(x)=\intrd e^{2\pi i \Phi\phas} \sigma\phas \hat{f}(\o)d\o.
\end{equation}
To avoid technicalities we
take $f\in\mathcal{S}(\R^d)$
or, more generally, $f\in
M^1$. If $\sigma\in L^\infty$
and the phase $\Phi$ is real,
the integral converges
absolutely and defines a
function in $L^\infty$.\par
Assume that the phase
function $\Phi(x,\eta)$
fulfills the following
properties:\\
(i) $\Phi\in \cC^{\infty}(\rdd)$;\\
(ii) for $z=\phas$,
\begin{equation}\label{phasedecay}
|\partial_z^\a \Phi(z)|\leq
C_\a,\quad |\a|\geq
2;\end{equation}
 (iii)
 there exists $\delta>0$ such
that
\begin{equation}\label{detcond}
    |\det\,\partial^2_{x,\eta} \Phi(x,\o)|\geq \delta.
\end{equation}
\par
If we set
\begin{equation}\label{cantra} \left\{
                 \begin{array}{l}
                 y=\nabla_{\eta}\Phi(x,\eta)
                 \\
                \xi=\nabla_{x}\Phi(x,\eta), \rule{0mm}{0.55cm}
                 \end{array}
                 \right.
\end{equation}
and solve with respect to $(x,\xi)$, we
obtain a
 mapping $\chi$,
defined by $(x,\xi)=\chi(y,\o)$, which
is a smooth bilipschitz canonical
transformation. This means that\\
-- $\chi$ is a  smooth diffeomorphism on $\rdd$;\\
-- both $\chi$ and $\chi^{-1}$ are Lipschitz continuous;\\
-- $\chi$ preserves the symplectic
form, i.e.,
$$dx\wedge d\xi=dy\wedge d\eta.$$
Indeed, under the above assumptions,
the global inversion function theorem
(see e.g. \cite{krantz}) allows us to
solve the first equation in
\eqref{cantra} with respect to $x$, and
substituting in the second equation
yields the smooth map $\chi$. The
bounds on the derivatives of $\chi$,
which give the Lipschitz continuity,
follow from the expression for the
derivatives of an inverse function
combined with the bounds in (ii) and
(iii). The symplectic nature of the map
$\chi$ is classical, see e.g.
\cite{caratheodory}. Similarly, solving
the second equation in \eqref{cantra}
with respect to $\eta$ one obtains the
function $\chi^{-1}$ with the desired
properties.\par In this section we
prove an almost diagonalization result
for FIOs as above, with respect to a
Gabor frame. Here we consider the case
of regular symbols. In Section
\ref{modsection} we will study the case
of symbols in modulation spaces.\par
Precisely, for a given $N\in\bN$, we
consider symbols $\sigma$ on $\rdd$
satisfying, for $z=\phas$,
\begin{equation}\label{decaysymbol}
    |\partial_z^\a \sigma(z)|\leq
    C_{\a},\quad\mbox{a.e.}\,\z\in\rdd,\ |\a|\leq 2
    N,
\end{equation}
here $\partial_z^\a$ denotes
distributional derivatives.\par

 Our goal is to study the decay properties of the matrix of the
 FIO $T$ with respect to a Gabor frame. For simplicity,
  we
 consider
 a normalized tight frame $\cG(g,\a,\b)$, with $g\in \cS(\rd).$
\begin{theorem}\label{T1}
Consider a phase function
satisfying {\rm (i)} and {\rm
(ii)} and a symbol satisfying
\eqref{decaysymbol}. There
exists $C_N>0$ such that
\begin{equation}\label{ET1}
    |\la T g_{m,n},\, g_{m',n'}\ra|\leq C_N {\la \nabla_z
    \Phi(m',n)-(n',m)\ra^{-2N}
    }.
\end{equation}
\end{theorem}
\begin{proof}
Recall that the time-frequency shifts
interchange under the action of the
\ft\,: $(T_x f)^{\wedge}=M_{-x}\hat{f}$
and $(M_{\o}
f)^{\wedge}=T_{\o}\hat{f}$, besides
they fulfill the commutation relations
$ T_x M_\o=e^{-2\pi i x\o} M_\o T_x.$
Using this properties, we can write
\begin{align*}
  \la T g_{m,n},&  g_{m',n'}\ra \\
   =& \intrd T g_{m,n}(x)\overline{M_{n'}T_{m'}g(x)}\,dx \\
   =& \intrd\intrd e^{2\pi i \Phi\phas}\sigma\phas T_{n}M_{-m}\hat{g}(\o)M_{-n'}T_{m'}\bar{g}(x)\,dx
   d\o\\
   =&\intrd\intrd M_{(0,-m)}T_{(0,-n)}\left(e^{2\pi i \Phi\phas}\sigma\phas\right) \hat{g}(\o)M_{-n'}T_{m'}\bar{g}(x)\,dx
   d\o\\
   =&\intrd\intrd T_{(-m',0)} M_{-(n',0)}M_{(0,-m)}T_{(0,-n)}\left(e^{2\pi i \Phi\phas}\sigma\phas\right) \bar{g}(x)\hat{g}(\o)\,dx
   d\o\\
    =&\intrd\intrd e^{2\pi i [\Phi(x+m',\o+n)-(n',m)\cdot (x+m',\o)]}\sigma(x+m',\o+n) \bar{g}(x)\hat{g}(\o)\,dx
   d\o
\end{align*}

Since $\Phi$ is smooth, we expand
$\Phi\phas$ into a Taylor series around
$(m',n)$ and obtain
$$\Phi(x+m',\o+n)=\Phi(m',n)+\nabla_z\Phi(m',n)\cdot (x,\o)+\Phi_{2,(m',n)}\phas$$
where the remainder is given by
 \begin{equation}
  \label{eq:c11}
\Phi_{2,(m',n)}\phas=2\sum_{|\a|=2}\int_0^1(1-t)\partial^\a
\Phi((m',n)+t\phas)\,dt\frac{\phas^\a}{\a!}.
\end{equation}
Whence, we can write
\begin{align*}
  |\la T g_{m,n},& \, g_{m',n'}\ra| \\
    =&\left|\intrd\intrd e^{2\pi i \{[\nabla_z\Phi(m',n)-(n',m)]\cdot \phas\}}e^{2\pi
    i\Phi_{2,(m',n)}\phas} \sigma(x+m',\o+n)
    \bar{g}(x)\hat{g}(\o)\,dx   d\o\right|
\end{align*}
For $N\in\bN$, using the identity:
$$(1-\Delta_z)^N e^{2\pi i \{[\nabla_z\Phi(m',n)-(n',m)]\cdot \phas\}}=\la
    2\pi(\nabla_z\Phi(m',n)-(n',m))\ra^{2N} e^{2\pi i \{[\nabla_z\Phi(m',n)-(n',m)]\cdot \phas\}},
$$
we integrate by parts and obtain
\begin{align*}
  |\la T g_{m,n},& \, g_{m',n'}\ra| \\
    =&\frac1{\la
    2\pi(\nabla_z\Phi(m',n)-(n',m))\ra^{2N}}\left|\intrd\intrd e^{2\pi i \{[\nabla_z\Phi(m',n)-(n',m)]\cdot
    \phas\}}\right.\\
    &\quad\quad
\left. \,\times\,(1-\Delta_z)^N\left[
    e^{2\pi
    i\Phi_{2,(m',n)}\phas} \sigma(x+m',\o+n)
    \bar{g}(x)\hat{g}(\o)\right]\,dx   d\o\right|.
\end{align*}
By means of Leibniz's formula the
factor  $$(1-\Delta_z)^N\left[
    e^{2\pi
    i\Phi_{2,(m',n)}\phas} \sigma(x+m',\o+n)
    \bar{g}(x)\hat{g}(\o)\right]$$
can be expressed as
$$e^{2\pi
    i\Phi_{2,(m',n)}(z)}\sum_{|\a|+|\b|+|\gamma|\leq 2N}
     C_{\a,\b\gamma}p( \partial^{|\a|} \Phi_{2,(m',n)})(z)
    (\partial_z^\b
    \sigma)(z+(m',n))\partial_z^\gamma(\bar{g}
    \otimes\hat{g})(z),$$
    where $p(\partial^{|\a|} \Phi_{2,(m',n)})(z)$ is a   polynomial  made  of derivatives of
    $\Phi_{2,(m',n)}$ of order at most $|\a|$.\par
    As a consequence of (ii) we have
$\partial_z^{\alpha}
\Phi_{2,(m',n)}(z)=O(\langle
z\rangle^{2})$, which
combined with the assumption
\eqref{decaysymbol} and the
hypothesis $g\in
\mathcal{S}(\R^d)$ yields the
desired estimate.
\end{proof}
\begin{remark}\rm
More generally one can consider symbols
satisfying estimates of the form
\begin{equation}\label{decaysymbol2}
    |\partial_z^\a \sigma(z)|\leq
    C_{\a}\mu(z),\quad\mbox{a.e.}\,z\in\rdd,\ |\a|\leq 2
    N,
\end{equation}
with $\mu\in\cM_v$ and also more
general windows $g$. Indeed, by arguing
as above and using
$$|\partial_z^\b
    \sigma(z+(m',n))|\leq C_{\b}\mu(z+(m',n))\leq  C' C_{\b}
    v(z) \mu(m',n),$$
one deduces the decay estimates
\begin{equation*}
  |\la T g_{m,n},  g_{m',n'}\ra| \leq C_N \frac{\mu(m',n)}{{\la
    \nabla_z\Phi(m',n)-(n',m)\ra^{2N}}},
\end{equation*}
 provided the integral
 \begin{equation}\label{convint}
\intrdd \left|
\sum_{|\a|+|\b|+|\gamma|\leq 2N}
C_{\a,\b,\gamma}p(
\partial^{|\alpha|} \Phi_{2,(m',n)})(z)
    v(z)\partial_z^\gamma(\bar{g}\otimes\hat{g})(z)
    \right|dz,
    \end{equation}
    converges. This is guaranteed
    if, e.g.,
    $\langle z\rangle ^{2N}
    v(z)\partial_z^\gamma(\bar{g}\otimes\hat{g})(z)\in
    L^1$.
\end{remark}
We now assume the additional hypothesis
(iii) on the phase, and rewrite
\eqref{ET1} in a form convenient for
the applications to the continuity of
FIOs in the next section. We need the
following lemma.
\begin{lemma}
Consider a phase function
$\Phi$ satisfying {\rm (i)},
{\rm (ii)}, and {\rm (iii)}.
Then
\begin{equation}
|\nabla_x\Phi(m',n)-n'|+|\nabla_\eta\Phi(m',n)-m|
\gtrsim|x(m,n)-m'|+|\xi(m,n)-n'|,
\end{equation}
where $(y,\eta)\longmapsto (x,\xi)$ is
the canonical transformation generated
by $\Phi$.
\end{lemma}
\begin{proof}
It suffices to prove the following
inequalities:
\begin{equation}\label{1f}
|\nabla_\eta\Phi(m',n)-m|\gtrsim
|x(m,n)-m'|,
\end{equation}
\begin{equation}\label{2f}
|\nabla_x\Phi(m',n)-n'|\geq
|\xi(m,n)-n'|-C|\nabla_\eta\Phi(m',n)-m|.
\end{equation}
We observe that, by \eqref{cantra}, we
have
\begin{equation}\label{3f}
y=\nabla_\eta
\Phi(x(y,\eta),\eta)\quad
\forall
(y,\eta)\in\mathbb{R}^{2d}
\end{equation}
and
\begin{equation}
\label{4f} \nabla_x\Phi
(x,\eta)=\xi(\nabla_\eta\Phi
(x,\eta),\eta)\quad\forall
(x,\eta)\in\mathbb{R}^{2d}.
\end{equation}
Hence, we have $m=
\nabla_\eta(x(m,n),n)$, so that
\begin{align*}
|\nabla_\eta
\Phi(m',n)-m|&=|\nabla_\eta\Phi(m',n)-\nabla_\eta
\Phi(x(m,n),n)|\\
&\gtrsim |x(m,n)-m'|,
\end{align*}
wehere the last inequality follows from
the fact that, for every fixed $\eta$,
the map $x\longmapsto
\nabla\Phi_\eta(x,\eta)$ has a
Lipschitz inverse, with Lipschitz
constant uniform with respect to
$\eta$.\par
 This proves \eqref{1f}.
\par
In order to prove \eqref{2f} we observe
that, in view of \eqref{4f}, it turns
out
\begin{align*}
\nabla_x\Phi(m',n)-n'&=\xi(\nabla_\eta\Phi(m',n),n)-n'\\
&=\xi(m+\nabla_\eta\Phi(m',n)-m,n)-n'\\
&=\xi(m,n)-n'+O(\nabla_\eta\Phi(m',n)-m)).
\end{align*}
where the last inequality follows from
the Taylor formula for the function
$y\longmapsto\xi(y,n)$, taking into
account that the function $\xi$ has
bounded derivatives. \par This proves
\eqref{2f}.
 \end{proof}\par
 Combining the
previous lemma with \eqref{ET1} we
obtain the following result.
\begin{theorem}\label{well-pose}
Consider a phase function
$\Phi$ satisfying {\rm (i)},
{\rm (ii)}, and {\rm (iii)},
and a symbol satisfying
\eqref{decaysymbol}. Let
$g\in\mathcal{S}(\R^d)$.
There exists a constant
$C_N>0$ such that
\begin{equation}\label{5f}
    |\la T g_{m,n}, g_{m',n'}\ra|\leq C_N {\la \chi(m,n)
    -(m',n')\ra^{-2N}
    },
\end{equation}
where $\chi$ is the canonical
transformation generated by $\Phi$.
\end{theorem}
This result shows that the
matrix representation of a
FIO with respect a Gabor
frame is well-organized,
similarly to the results
recently obtained by
\cite{candes,Guo-labate} in
terms of shearlets and
curvelets frames. More
precisely, if $\sigma\in
S^0_{0,0}$, namely if
\eqref{decaysymbol} is
satisfied for every
$N\in\mathbb{N}$, then the
Gabor matrix of $T$ is highly
concentrated along the graph
of $\chi$.

 \section{Continuity of FIOs on $M^p_\mu$}
 In this section we study the continuity of FIOs on
 the modulation spaces $M^p_\mu$ associated with a
 weight function $\mu\in\cM_{v_s}$,
  $s\geq0$.
 We need the
   following preliminary lemma.
   \begin{lemma}\label{schur}
   Consider a lattice $\Lambda$ and an
   operator $K$ defined on sequences as
   \[
   (Kc)_{\lambda}=\sum_{\nu\in\Lambda}
   K_{\lambda,\nu}c_\nu,
   \]
   where
   \[
   \sup_{\nu\in\Lambda}\sum_{\lambda\in\Lambda}
   |K_{\lambda,\nu}|<\infty,\quad \sup_{\lambda\in\Lambda}
   \sum_{\nu\in\Lambda}
   |K_{\lambda,\nu}|<\infty.
   \]
   Then $K$ is continuous on
   $l^p(\Lambda)$ for every $1\leq
   p\leq\infty$ and moreover maps
   the space $c_0(\Lambda)$ of sequences vanishing at infinity into
   itself.
   \end{lemma}
   \begin{proof}
The first part is the
classical Schur's test (see
e.g. \cite[Lemma
6.2.1]{book}). The second
part follows in this way.
Since we know that $K$ is
continuous on $l^\infty$ and
the space of eventually zero
sequences is dense in $c_0$,
it suffices to verify that
$K$ maps every eventually
zero sequence in $c_0$. This
follows from the fact that
any eventually zero sequence
belongs to $l^1$ and
therefore, since $K$ is
continuous on $l^1$, is
mapped in $l^1\hookrightarrow
c_0$.\end{proof}\\
 We can now state our result.
\begin{theorem}\label{contmp}
Consider a phase function
satisfying {\rm (i)}, {\rm
(ii)}, and {\rm (iii)}, and a
symbol satisfying
\eqref{decaysymbol}. Let
$0\leq s<2N-2d$, and
$\mu\in\cM_{v_s}$. For every
$1\leq p<\infty$, $T$ extends
to a continuous operator from
${M}^p_{\mu\circ\chi}$ into
${M}^p_{\mu}$, and for
$p=\infty$ it extends to a
continuous operator from
$\tilde{{M}}^\infty_{\mu\circ\chi}$
into
$\tilde{{M}}^\infty_{\mu}$.
\end{theorem}
Recall that
$\tilde{{M}}^\infty_{\mu}$ is
the closure of
$\mathcal{S}(\R^d)$ in
${{M}}^\infty_{\mu}$.
Moreover, observe that
$\mu\circ\chi\in\cM_{v_s}$.
Indeed,
 $v_s\circ\chi\asymp v_s$,  due to the bilipschitz property of $\chi$.

\begin{proof}
We first prove that
\[
\|Tf\|_{M^{p}_{\mu\circ\chi}}\leq
C\|f\|_{M^p_\mu},
\]
for every $f\in
\mathcal{S}(\R^d)$. This
proves the theorem in the
case $p<\infty$, since
$\mathcal{S}(\R^d)$ is dense
in $M^p_\mu$.\par
 We see at once that, since $\sigma\in L^\infty$, $T$
defines a bounded operator
from $M^1$ into
$L^\infty\hookrightarrow
M^\infty$. Hence, for all
$f\in \mathcal{S}(\R^d)$, we
have $Tf\in M^\infty$ and
Theorem \ref{teomod} shows
that
$\|f\|_{M^p_{\mu\circ\chi}}\gtrsim
\|C_g(f)\|_{l^p_{\mu\circ\chi}}$
and $\|Tf\|_{M^p_\mu}\lesssim
\|C_g(Tf)\|_{l^p_\mu}$.  On
the other hand, the expansion
\eqref{parsevalframe} holds
for $f$ with convergence in
$M^1$. Therefore
\[
Tf=\sum_{m,n}\langle f,
g_{m,n}\rangle Tg_{m,n}\]
with convergence in
$M^\infty$. Hence,
\[
C_g(Tf)_{m',n'}=\langle
Tf,g_{m',n'}\rangle=\sum_{m,n}\langle T
g_{m,n},g_{m',n'}\rangle \langle
f,g_{m,n}\rangle=\sum_{m,n}\langle T
g_{m,n},g_{m',n'}\rangle C_g(f)_{m,n}.
\]
Therefore we are reduced to proving
that the matrix operator
\begin{equation}\label{matop}
\{c_{m,n}\}\longmapsto
\sum_{m,n\in\mathbb{Z}^d} \langle T
g_{m,n},g_{m',n'}\rangle c_{m,n}
\end{equation}
is bounded from
$l^{p}_{\widetilde{\mu\circ\chi}}$
into $l^p_{\tilde{\mu}}$.
This follows from Schur's
test (Lemma \ref{schur}) if
we prove that, upon setting
\[
K_{m',n',m,n}=\langle T
g_{m,n},g_{m',n'}\rangle
\mu(m',n')/\mu(\chi(m,n)),
\]
we have
\begin{equation}\label{6f}
K_{m',n',m,n}\in
l^\infty_{m,n}l^1_{m',n'},
\end{equation}
and
\begin{equation}\label{7f}
K_{m',n',m,n}\in
l^\infty_{m',n'}l^1_{m,n}.
\end{equation}
In view of \eqref{5f} we have
\begin{equation}\label{8f}
|K_{m',n',m,n}|\lesssim
\langle\chi(m,n)-(m',n')\rangle^{-2N+s}
\frac{\mu(m',n')}{\langle\chi(m,n)-(m',n')\rangle^s
\mu(\chi(m,n))}.
\end{equation}
Now, the last quotient in
\eqref{8f} is bounded because
$\mu$ is ${v_s}$-moderate, so
we deduce $\eqref{6f}$.\par
Finally, since $\chi$ is a
bilipschitz function we have
\begin{equation}\label{9f}
|\chi(m,n)-(m',n')|\asymp|(m,n)-\chi^{-1}(m',n')|
\end{equation}
so that \eqref{7f} follows as well.\par
The case $p=\infty$ follows analogously
by using Theorem \ref{teomod2} (with
$p=q=\infty$), and the last part of the
statement of Lemma \ref{schur}.
\end{proof}
\begin{remark}\rm
Theorem \ref{contmp} with
$v\equiv1$ gives, in
particular, continuity on the
unweighted modulation spaces
$M^p$. If moreover $p=2$, we
recapture the classical
$L^2$-continuity result by
Asada and Fujiwara
\cite{asada-fuji}.\\
Also, Theorem \ref{contmp}
applies to $\mu=v_t$, with
$|t|\leq s$. In that case we
obtain continuity on
$M_{v_t}$, because
$v_t\circ\chi\asymp v_t$.
\end{remark}

 \section{Continuity of FIOs on $\mpq$}
In this section we study the continuity
of FIOs on modulation spaces $M^{p,q}$
possibly with $p\not=q$. As shown in
Section \ref{metaplectic}, under the
assumptions of Theorem \ref{contmp}
such operators may fail to be bounded
when $p\not=q$. The counterexample is
given by the phase
$\Phi(x,\eta)=x\eta+|x|^2/2$, and
symbol $\sigma=1$, which does not yield
a bounded operator on $M^{p,q}$, except
for the case $p=q$. Here the
obstruction is essentially due to the
fact that the map
$x\longmapsto\nabla_x\Phi(x,\eta)$ has
unbounded range. Indeed we will show,
for general phases, that if such a map
has range of finite diameter, uniformly
with respect to $\eta$, then the
corresponding operator is bounded on
all $M^{p,q}$. To this end we need the
following result.
\begin{proposition}\label{proschur}
Consider an operator defined on
sequences on the lattice
$\Lambda=\alpha\mathbb{Z}^d\times\beta\mathbb{Z}^d$
by
\[
(Kc)_{m',n'}=\sum_{m,n}
{K_{m',n',m,n}}c_{m,n}.
\]
(i) If $K\in
l^\infty_{n}l^1_{n'}l^\infty_{m'}l^1_{m}$,
$K$ is continuous on
$l^1_{n}l^\infty_{m}$.\\
(ii) If $K\in
l^\infty_{n'}l^1_{n}l^\infty_{m}l^1_{m'}$,
$K$ is continuous on
$l^\infty_{n}l^1_{m}$.\\
(iii) If $K\in
l^\infty_{n}l^1_{n'}l^\infty_{m'}l^1_{m}\cap
l^\infty_{n'}l^1_{n}l^\infty_{m}l^1_{m'}$,
and moreover $K\in
l^\infty_{m',n'}l^1_{m,n}\cap
l^\infty_{m,n}l^1_{m',n'}$,
the operator $K$ is
continuous on
$l^{p,q}=l^q_{n}l^p_{m}$ for
every $1\leq p,q\leq\infty$.\\
(iv) Assume the hypotheses in (iii).
Then $K$ is continuous on all
$\tilde{l}^{p,q}$, $1\leq
p,q\leq\infty$.
\end{proposition}
Recall that $\tilde{l}^{p,q}$
is the closure of the space
of eventually zero sequences
in ${l}^{p,q}$.
\begin{proof}
$(i)$ We have
\begin{align*}
\|Kc\|_{l^1_{n'}l^\infty_{m'}}&\leq\sum_{n'}
\sup_{m'}\sum_{m,n}
|{K_{m',n',m,n}}||c_{m,n}|\\
&\leq\sum_n\left(\sum_{n'}
\sup_{m'}\sum_m
|{K_{m',n',m,n}}|\right)\sup_{m}|c_{m,n}|\\
&\leq
\|K\|_{l^\infty_{n}l^1_{n'}l^\infty_{m'}l^1_{m}}
\|c\|_{l^1l^\infty}.
\end{align*}
$(ii)$ It turns out
\begin{align*}
\|Kc\|_{l^\infty_{n'}l^1_{m'}}&\leq
\sup_{n'}\sum_{m'}\left(\sum_{m,n}|{K_{m',n',m,n}}||c_{m,n}|\right)\\
&\leq\sup_{n'}\sum_n\left(\sup_{m}\sum_{m'}
|{K_{m',n',m,n}}|
\sum_m
|c_{m,n}|\right)\\
&\leq
\|K\|_{l^\infty_{n'}l^1_{n}l^\infty_{m}l^1_{m'}}\|c\|_{l^\infty
l^1}.
\end{align*}
$(iii)$ Since the statement
holds for $p=q$ by the
classical Schur's test, and
for $(p,q)=(1,\infty)$ and
$(p,q)=(\infty,1)$ by the
items $(i)$ and $(ii)$, it
follows by complex
interpolation (see (3) on
page 128 and (15) on page 134
of \cite{triebel}) that it
holds for all $(p,q)$, except
possibly in the cases
$q=\infty$, $1<p<\infty$. For
these cases we argue by
duality as follows.\par In
order to prove the continuity
of $K$ on $l^\infty l^p$, it
suffices to verify that for
any sequences $c=(c_{m,n})\in
l^\infty_n l^p_m$ and
$d=(d_{m',n'})\in
l^1_{n'}l^{p'}_{m'}$, with
$d_{m',n'}\geq0$, we have
\begin{equation}\label{ddua}
\sum_{m',n'}|(Kc)_{m',n'}|d_{m',n'}\lesssim
\|c\|_{l^\infty_n
l^p_m}\|d\|_{l^1_nl^{p'}_m}.
\end{equation}
Now
\begin{align*}
\sum_{m',n'}|(Kc)_{m',n'}|d_{m',n'}&\leq
\sum_{m',n'}\sum_{m,n}|K_{m',n',m,n}||c_{m,n}|d_{m',n'}\\
&=\sum_{m,n}\left(\sum_{m',n'}|K_{m',n',m,n}|
d_{m',n'}\right)|c_{m,n}|\\
&\leq \|c\|_{l^\infty
l^p}\|\tilde{K}d\|_{l^1l^{p'}}.
\end{align*}
where $\tilde{K}$ is the
operator with matrix kernel
$\tilde{K}_{m,n,m',n'}=|K_{m',n',m,n}|$.
Since it satisfies the same
assumptions as $K$, it is
continuous on $l^1l^{p'}$,
which gives \eqref{ddua}.\\
$(iv)$ Since $K$ is
continuous on $l^{p,q}$ and
by the definition of
$\tilde{l}^{p,q}$, it
suffices to verify that $K$
maps every eventually zero
sequence in
$\tilde{l}^{p,q}$. This
follows from the fact that
$K$ maps every eventually
zero sequence in
$l^1\hookrightarrow\tilde{l}^{p,q}$,
because $K$ is bounded on
$l^1$.
\end{proof}

\begin{theorem}\label{cont}
Consider a phase function
$\Phi$ satisfying {\rm (i)},
{\rm (ii)}, and {\rm (iii)},
and a symbol satisfying
\eqref{decaysymbol}, with
$N>d$. Suppose, in addition,
that
\begin{equation}\label{fase}
\sup_{x,x',\eta}\left|
\nabla_x\Phi(x,\eta)-\nabla_x\Phi(x',\eta)\right|<\infty.
\end{equation}
Then the corresponding Fourier integral
operator $T$ extends to a bounded
operator on $M^{p,q}$ for every $1\leq
p,q<\infty$ and on $\tilde{M}^{p,q}$ if
$p=\infty$ or $q=\infty$.
\begin{proof}
By arguing as in the proof of Theorem
\ref{contmp}, it suffices to prove the
continuity on $l^{p,q}=l^{q}_n l^p_m$
if $p<\infty$ and $q<\infty$, or
$\tilde{l}^{p,q}$ if $p=\infty$ or
$q=\infty$, of the operator
\[
\{c_{m,n}\}\longmapsto\sum_{m,n\in\mathbb{Z}^{d}}
T_{m',n',m,n} c_{m,n},
\]
where
\[
T_{m',n',m,n}=\langle
Tg_{m,n},g_{m',n'}\rangle.
\]
By applying \ref{proschur},
it suffices to verify that
\begin{equation}\label{uno}
\{T_{m',n',m,n}\}\in l^\infty_n
l^1_{n'} l^\infty _{m'} l^1_m,
\end{equation}
\begin{equation}\label{due}
\{T_{m',n',m,n}\}\in l^\infty_{n'}
l^1_{n} l^\infty _{m} l^1_{m'},
\end{equation}
because we already see from
\eqref{5f} and \eqref{9f}
that $\{T_{m',n',m,n}\}\in
l^\infty_{m,n}
l^1_{m',n'}\cap
l^\infty_{m',n'}
l^1_{m,n}$.\par Let us now
prove \eqref{uno}. It follows
from \eqref{ET1} and
\eqref{2f} that
\begin{align*}
|T_{m',n',m,n}|&\lesssim
\left(1+\left|\nabla_x\Phi(m',n)-n'\right|^2
+\left|\nabla_\eta\Phi(m',n)-m\right|^2\right)^{-N}\\
&\lesssim
\left(1+\left|\xi(m,n)-n'\right|^2
+\left|\nabla_\eta\Phi(m',n)-m\right|^2\right)^{-N}\\
&\lesssim
\left(1+\left|\xi(m,n)-n'\right|\right)^{-N}
\left(1+\left|\nabla_\eta\Phi(m',n)-m\right|\right)^{-N}.
\end{align*}
By \eqref{cantra} we have
\[
\xi(y,\eta)=\nabla_x\Phi(x(y,\eta),\eta),\quad
\forall(y,\eta)\in\mathbb{R}^{2d},
\]
so that the hypothesis
\eqref{fase} yields
\[
\xi(m,n)=\xi(0,n)+O(1).
\]
Hence \eqref{uno} follows.\par We now
prove \eqref{due}. As above, it follows
from \eqref{5f} and \eqref{1f} that
\begin{align*}
|T_{m',n',m,n}|&\lesssim
\left(1+\left|\nabla_x\Phi(m',n)-n'\right|^2
+\left|\nabla_\eta\Phi(m',n)-m\right|^2\right)^{-N}\\
&\lesssim
\left(1+\left|\nabla_x\Phi(m',n)-n'\right|^2
+\left|x(m,n)-m'\right|^2\right)^{-N}\\
&\lesssim
\left(1+\left|\nabla_x\Phi(m',n)-n'\right|\right)^{-N}
\left(1+\left|x(m,n)-m'\right|\right)^{-N}
\end{align*}
By \eqref{fase} we have
\[
\nabla_x\Phi(m',n)=\nabla_x\Phi(0,n)+O(1),
\]
so that
\begin{align}
1+\left|\nabla_x\Phi(m',n)-n'\right|&\gtrsim
1+\left|\nabla_x\Phi(0,n)-n'\right|\nonumber\\
&\gtrsim
1+|n-\psi(n')|,\label{nuo}
\end{align}
where $\psi$ is the inverse
function of the  bilipschitz
function
$\eta\longmapsto\nabla_x\Phi(0,\eta)$.
Therefore we obtain
\eqref{due}.\par
 This concludes the proof.
\end{proof}
\end{theorem}
\begin{example}\label{osser3}\rm
Theorem \ref{cont} applies,
in particular, to phases of
the type
\[
\Phi(x,\eta)=x\eta+a(x,\eta)\
\ {\rm where}\
|\partial^\alpha_x\partial^\beta_\eta
a(x,\eta)|\leq
C_{\alpha,\beta}\ {\rm for\,
every \,}\
2|\alpha|+|\beta|\geq2.
\]
In the special case when
$a(x,\eta)=a(\eta)$ is
independent of $x$ and the symbol $\sigma\equiv 1$, the FIO reduces to a Fourier multiplier
$$Tf(x)= \intrd e^{2\pi i x\eta} e^{2\pi i a(\eta)}\hat{f}(\eta)\,d\eta
$$
and we reobtain the
result of \cite[Theorem 5]{benyi} on
the continuity of $T$ on all $M^{p,q}$,
$1\leq p,q\leq\infty$.
\end{example}

\section{Modulation spaces as symbol
classes}\label{modsection}

In what follows we shall rephrase the
quantity $|\la T g_{m,n},
g_{m',n'}\ra|$ in terms of the STFT of
the symbol $\sigma$, without assuming
the existence of derivatives of
$\sigma$. This will be applied  to prove the continuity
 of FIOs with symbols in $M^{\infty,1}$ on modulation spaces $M^p$.
\par The same arguments as in Theorem \ref{T1} yield
the equality
$$\la Tg_{m,n},g_{m',n'}\ra=e^{2\pi i mn}\intrd e^{2\pi i \Phi\phas} \sigma\phas M_{(-n',-m)}T_{(m',n)}(\bar{g}\otimes\hat{g})\phas\,dx d\o.$$
Expanding the phase $\Phi$  into a
Taylor series around $(m',n)$ we obtain

$$\Phi(x,\o)=\Phi(m',n)+\nabla_z\Phi(m',n)\cdot (x-m',\o-n)+T_{(m',n)}\Phi_{2,(m',n)}\phas$$
where the remainder $\Phi_{2,(m',n)}$
is given by
 \eqref{eq:c11}.\par
Inserting this expansion in the
integrals above, we can write
\begin{align}\label{F1}
\la Tg_{m,n},g_{m',n'}\ra&=e^{2\pi i (mn+\Phi(m',n)
-\nabla_z\Phi(m',n)\cdot(m',n))}\intrd e^{2\pi i \nabla_z\Phi(m',n)\phas}\sigma\phas\\
&\quad\quad\quad\times\quad
M_{(-n',-m)}T_{(m',n)}e^{2\pi i
\Phi_{2,(m',n)}\phas}(\bar{g}\otimes\hat{g})\phas\,dx
d\o.\nonumber
\end{align}
Defining
\begin{equation}\label{window}
\Psi_{(m',n)}\phas:=e^{2\pi
i\Phi_{2,(m',n)}\phas}(\bar{g}\otimes\hat{g})\phas,
\end{equation}
and computing the modulus of the
left-hand side of \eqref{F1}, we are
led to

\begin{equation}\label{symbstft}|\la T g_{m,n},  g_{m',n'}\ra| = |V_{\Psi_{(m',n)}} \sigma((m',n),(n'-\nabla_x\Phi(m',n),m-\nabla_\o \Phi(m',n)))|
\end{equation}

Observe that the window $\Psi_{(m',n)}$
of the STFT above depends on the pair
$(m',n)$.

We now
study the continuity problem of $T$
when the symbol $\sigma$ is in the modulation space
$M^{\infty,1}$.
\begin{theorem}\label{contmp2}
Consider a phase function
satisfying {\rm (i)}, {\rm
(ii)}, and {\rm (iii)}, and a
symbol $\sigma\in
M^{\infty,1}$. For every
$1\leq p<\infty$, $T$ extends
to a continuous operator on
${M}^p$, and for $p=\infty$
it extends to a continuous
operator on
$\tilde{{M}}^\infty$.
\end{theorem}
By arguing as in the proof of
Theorem \ref{contmp}, it
suffices to prove the
continuity on $l^p$ if $1\leq
p<\infty$ and on
$\tilde{l}^\infty=c_0$, of
the operator \eqref{matop}.
In wiew of Schur's test
(Lemma \ref{schur}) and
\eqref{symbstft}, it suffices
to prove the following
result.
\begin{proposition}\label{T2}
Consider a phase function
$\Phi$ satisfying {\rm (i)}
and {\rm (ii)} and {\rm
(iii)} and a symbol
$\sigma\in M^{\infty,1}$. If
we set
\begin{equation}\label{zmn}
   z_{m,n,m',n'}:= ((m',n),(n'-\nabla_x\Phi(m',n),m-\nabla_\o \Phi(m',n))),\quad m,m'\in \a\zd,\, n,n'\in \b\zd,
\end{equation}
then,
\begin{eqnarray}
 \sup _{(m,n)\in\Lambda}\sum_{(m',n')\in\Lambda} |V_{\Psi_{(m',n)}}\sigma( z_{m,n,m',n'})|&\lesssim& \|\sigma\|_{M^{\infty,1}}.\label{EW1}\\
   \sup _{(m',n')\in\Lambda}\sum_{(m,n)\in\Lambda} |V_{\Psi_{(m',n)}}\sigma( z_{m,n,m',n'})|&\lesssim& \|\sigma\|_{M^{\infty,1}}.\label{EW2}
  \end{eqnarray}
\end{proposition}
We need the following lemma.
\begin{lemma}\label{changewind2}
Let $\Psi_0\in\cS(\rdd)$ with
$\|\Psi_0\|_{L^2}=1$ and
$\Psi_{(m',n)}$ be defined by
\eqref{window}, with
$(m',n)\in
\Lambda=\a\zd\times\b\zd$,
and $g\in\cS(\rd)$. Then,
\begin{equation}\label{newwind}
\int_{\bR^{4d}}\sup
_{(m',n)\in\Lambda}|V_{\Psi_{(m',n)}}\Psi_0(w)|\,dw<\infty.
\end{equation}
\end{lemma}
\begin{proof}[Proof of Lemma
\ref{changewind2}] We shall
show that
\begin{equation}\label{lw} |V_{\Psi_{(m',n)}}\Psi_0(w)|\leq C\la w\ra^{-(4d+1)},\quad \quad \forall (m',n)\in\Lambda.
\end{equation}
Using the switching property of the
STFT: $$ (V_fg)(x,\o)=e^{-2\pi i\o
x}\overline{(V_gf)(-x,-\o)},$$ we
observe that
$|V_{\Psi_{(m',n)}}\Psi_0(w_1,w_2)|=|V_{\Psi_0}\Psi_{(m',n)}|(-w_1,-w_2)$,
and by the even property of of the
weight $\la \cdot\ra$, relation
\eqref{lw} is equivalent to
\begin{equation}\label{l1} |V_{\Psi_0}\Psi_{(m',n)}(w)|\leq C\la w\ra^{-(4d+1)},\quad \quad \forall (m',n)\in\Lambda.
\end{equation}
Now, the mapping $V_{\Psi_0}$ is
continuous from $\cS(\rdd)$ to
$\cS(\bR^{4d})$ (see \cite[Chap.
11]{book}). This means that there
exists $M\in\bN$, $K>0$,  such that
\begin{align*}|V_{\Psi_0}\Psi_{(m',n)}(w)|&=|V_{\Psi_0}\Psi_{(m',n)}(w)|\frac{\la w\ra^{4d+1}}{\la w\ra^{4d+1}}\\
&\leq \|V_{\Psi_0}\Psi_{(m',n)}\la \cdot\ra^{4d+1}\|_{L^\infty(\bR^{4d})}\la w\ra^{-(4d+1)}\\
&\leq K \sum_{|\gamma|+|\delta|\leq
M}\|D^\gamma X^\delta
\Psi_{(m',n)}\|_{L^\infty(\rdd)}\la
w\ra^{-(4d+1)},
\end{align*}
for every $(m',n)\in\Lambda$.
We now claim that
$\Psi_{(m',n)}\in\cS(\rdd)$
uniformly with respect to
$(m',n)$. This is proved as
follows: the function
$e^{2\pi i
\Phi_{2,(m',n)}\phas}$ is in
$\cC^{\infty}(\rdd)$ and
possesses derivatives
dominated by powers
$\la\phas\ra^k$, $k\in\N$,
uniformly with respect to
$(m',n)$, due to
\eqref{phasedecay}; since
$(\bar{g}\otimes\hat{g})\in
\cS(\rdd)$, it follows that
$\Psi_{(m',n)}\in\cS(\rdd)$,
with semi-norms  uniformly
bounded:
$$p_{M}(\Psi_{(m',n)}):=\sum_{|\gamma|+|\delta|\leq M}\|D^\gamma X^\delta \Psi_{(m',n)}\|_{L^\infty}\leq C_{M},\quad\quad\forall (m',n)\in\Lambda.$$
Consequently,
$$ |V_{\Psi_0}\Psi_{(m',n)}(w)|\leq K  \sum_{|\gamma|+|\delta|\leq M}\|D^\gamma X^\delta \Psi_{(m',n)}\|_{L^\infty(\rdd)}\la w\ra^{-(4d+1)}\leq K C_{M}\la w\ra^{-(4d+1)},$$
for every $(m',n)\in\Lambda$, as
desired.
\end{proof}
\begin{proof}[Proof of Proposition \ref{T2}] We shall prove \eqref{EW1}. First,  Lemma \ref{changewind} for $g_1=\gamma=\Psi_0$, yields
$$|V_{\Psi_{(m',n)}}\sigma(z)|\leq (|V_{\Psi_0} \sigma |\ast|V_{\Psi_{(m',n)}}\Psi_0|)(z),\quad z\in\bR^{4d},
$$
so that
\begin{align*}
 \sum_{(m',n')\in\Lambda} &|V_{\Psi_{(m',n)}}\sigma( z_{m,n,m',n'})|\\
 &\leq\int_{\bR^{4d}}\sum_{(m',n')\in\Lambda}
 |V_{\Psi_0} \sigma (z_{m,n,m',n'}-w)|\,|V_{\Psi_{(m',n)}}\Psi_0(w)|\,dw\\
 &\leq \int_{\bR^{4d}}\sum_{(m',n')\in\Lambda}
 |V_{\Psi_0} \sigma (z_{m,n,m',n'}-w)|\,\sup _{(m',n)\in\Lambda}| V_{\Psi_{(m',n)}}\Psi_0(w)|\,dw\\
 & \leq \sup_{w\in\bR^{4d}}\sum_{(m',n')\in\Lambda}
 |V_{\Psi_0} \sigma (z_{m,n,m',n'}-w)|\,\int_{\bR^{4d}}\sup _{(m',n)\in\Lambda}| V_{\Psi_{(m',n)}}\Psi_0(w)|\,dw\\
 & \leq C \sup_{w\in\bR^{4d}}\sum_{(m',n')\in\Lambda}
 |V_{\Psi_0} \sigma (z_{m,n,m',n'}-w)|\,,
 \end{align*}
 where the last majorization is due to Lemma \ref{changewind2}. Since,
 \begin{equation*}
 \sum_{(m',n')\in\Lambda}
 |V_{\Psi_0} \sigma (z_{m,n,m',n'}-w)|\leq
 \sum_{(m',n')\in\Lambda}\sup_{u_1\in\rdd}|V_{\Psi_0}
  \sigma (u_1,\tilde{z}_{m,n,m',n',w_2})|,
 \end{equation*}
with
 $$\tilde{z}_{m,n,m',n',w_2}:=(n'-\nabla_x\Phi(m',n),
 m-\nabla_\o \Phi(m',n))-w_2,\quad
 w=(w_1,w_2)\in\bR^{4d},
 $$
we shall  prove that
 \begin{equation}\label{T2e}
 \sum_{(m',n')\in\Lambda}\sup_{u_1\in\rdd}|V_{\Psi_0}
 \sigma (u_1,\tilde{z}_{m,n,m',n',w_2})|\lesssim \|\sigma\|_{M^{\infty,1}},
 \end{equation}
 uniformly with respect to $(m,n)\in\Lambda$,
  $w_2\in \bR^{2d}$.
For every fixed $(m,n)$,  the
set $\cX=\cX_{m,n,w_2}$,
given by
$$ \cX_{m,n,w_2} = \{ \tilde{z}_{m,n,m',n',w_2};\,\,(m',n')\in \Lambda,\,w_2\in\rdd\},
$$
 is separated, uniformly with respect to $(m,n)$,
  $w_2$. Indeed, given $(m'_1,n'_1)\not=(m'_2,n'_2)$,
  if $m'_1\not=m'_2$,
 $$|\tilde{z}_{m,n,m'_1,n'_1,w_2}-\tilde{z}_{m,n,m'_2,
 n'_2,w_2}|\geq |\nabla_\o\Phi(m_1',n)-\nabla_\o \Phi(m_2',n)|\geq C |m_1'-m_2'|\geq \a C,
 $$
 uniformly with respect to $(m,n)$,  $w_2$,  because the mapping $x\longmapsto\nabla_\o\Phi(x,\o)$ has an inverse that is    Lipschitz continuous, thanks to \eqref{phasedecay} and \eqref{detcond}.
On the other hand, if $m'_1=m'_2$,
$$|\tilde{z}_{m,n,m'_1,n'_1,w_2}-\tilde{z}_{m,n,m'_2,n'_2,w_2}|
\geq  |n_1'-n_2'|\geq
\b,\quad \forall (m,n)\in
\Lambda, w_2\in\rdd.
$$
Hence, $\cX$ is separated
uniformly with respect to
$(m,n)$, $w_2$. Now, we apply
Proposition \ref{samp} (with
$p=q=1$) to the function
$$F(u_2)=\sup_{u_1\in\rdd}| V_{\Psi_0}
 \sigma(u_1,u_2)|,\quad u_2\in\bR^{2d},
$$
which is lower
semi-continuous, being
$V_{\Psi_0}\sigma$
continuous. We obtain
\begin{equation}\label{bo}\|F_{|_{\cX}}\|_{\ell^{1}}\leq C \|F\|_{W(L^\infty,L^{1})}=\|V_{\Psi_0} \sigma\|_{W(\cC,L^{\infty,1})}.
\end{equation}
If the symbol $\sigma$ is in
$M^{\infty,1}$, by Lemma
\ref{amalg} the STFT
$V_{\Psi_0} \sigma$ belongs
to the Wiener amalgam space
$W(\cF L^1,L^{\infty,1})$,
and
 \begin{equation*}
\| V_{\Psi_0}
\sigma\|_{W(\cC,L^{\infty,1})}\,
\lesssim  \| V_{\Psi_0}
\sigma\|_{W(\Fur
L^1,L^{\infty,1})}\lesssim \|\sigma\|_{
  M^{\infty,1}}\|\Psi_0\|_{M^{1}} .
\end{equation*}
The first inequality is due
to $\Fur
L^1\hookrightarrow\cC$ and
the inclusion relations
between Wiener amalgam
spaces. Combining this
inequality with \eqref{bo} we
obtain \eqref{T2e}, uniformly
with respect to $(m,n)$ and
$w_2$, that is
 \eqref{EW1}.\par   The estimate \eqref{EW2} is obtained by
 similar arguments.
\end{proof}
\begin{remark}\rm
We observe that the continuity on
$M^2=L^2$ of FIOs as above,
 with symbols in
$M^{\infty,1}$, was already proved in
\cite{boulkhemair} by other methods.
\end{remark}

 \section{The case of quadratic phases: metaplectic
 operators}\label{metaplectic}
In this section we briefly
discuss the particular case
of quadratic phases, namely
phases of the type
\begin{equation}\label{faseq}
\Phi(x,\eta)=\frac{1}{2}
Ax\cdot
x+Bx\cdot\eta+\frac{1}{2}C\eta\cdot\eta+\eta_0\cdot
x-x_0\cdot\eta,
\end{equation}
where
$x_0,\eta_0\in\mathbb{R}^d$,
$A,C$ are real symmetric
$d\times d$ matrices and $B$
is a real $d\times d$
nondegenerate matrix.\par It
is easy to see that, if we
take the symbol
$\sigma\equiv1$ and the phase
\eqref{faseq}, the
 corresponding FIO $T$ is (up to a constant
factor) a metaplectic
operator. This can be seen by
means of the easily verified
factorization
\begin{equation}\label{factorization}
 T=M_{\eta_0} U_A
D_B\Fur^{-1} U_C\Fur T_{x_0},
\end{equation}
where $U_A$ and $U_C$ are the
multiplication operators by
$e^{\pi i Ax\cdot x}$ and
$e^{\pi i C\eta\cdot \eta}$
respectively, and $D_B$ is
the dilation operator
$f\mapsto f(B\cdot)$. Each of
the factors is (up to a
constant factor) a
metaplectic operator (see
e.g. the proof of
\cite[Theorem
18.5.9]{hormander}), so $T$
is.\par The corresponding
canonical map, defined by
\eqref{cantra}, is now an
affine symplectic map. For
the benefit of the reader,
some important special cases
are detailed in the table
below. \vskip0.5truecm
\par
\begin{center}
\def\V{\rule{0pt}{2.5ex}}
\begin{tabular}{|l*{3}{|c}|}
\hline
 operator \V& phase $\Phi(x,\eta)$ & canonical transformation\\
\hline
 \ \ $T_{x_0}$\V& $(x-x_0)\cdot\eta$ & $\chi(y,\eta)=(y+x_0,\eta)$ \\
\hline
 \ \ $M_{\eta_0}$\V& $(\eta+\eta_0)\cdot x$ & $\chi(y,\eta)=(y,\eta+\eta_0)$ \\
\hline
\ \ $D_B$\V& $Bx\cdot\eta$ & $\chi(y,\eta)=(B^{-1}y,^t\!\!B\eta)$ \\
\hline
\ \ $U_A$\V& $x\cdot\eta+\frac{1}{2}Ax\cdot x$ & $\chi(y,\eta)=(y,\eta+Ax)$ \\
\hline
\end{tabular}
\end{center}
\vskip0.5truecm However one should
observe that there are metaplectic
operators, as the Fourier transform,
which cannot be expressed as FIOs of
the type \eqref{fio}.\par

 Metaplectic operators are
known to be bounded on
$M^p_{v_s}$, see e.g.
\cite[Proposition
12.1.3]{book}. This also
follows from Theorem
\ref{contmp}. Indeed, since
$\chi$ is a bilipschitz
function, we have
$v_s\circ\chi\asymp v_s$.

\par Also, Theorem \ref{cont} applies to
quadratic phases whose affine
symplectic map $\chi$ is (up
to translations on the phase
space) defined by an
upper-triangular matrix,
which happens precisely when
$A=0$. Indeed, we obtain the
map $\chi$ by solving
\begin{equation*} \left\{
                 \begin{array}{l}
                 y=Bx+C\eta-x_0
                 \\
                \xi=Ax+B\eta+\eta_0. \rule{0mm}{0.55cm}
                 \end{array}
                 \right.
\end{equation*}
The phase condition
\eqref{detcond} here becomes
\begin{equation*}
    \det\,B \not=0,
\end{equation*}
so  that  $B$ is an
invertible matrix and
$x=B^{-1}y-B^{-1}
C\eta+B^{-1}x_0$.
 Whence, the mapping  $\chi :(y,\eta)\longmapsto (x,\xi)$ is given by
$$\bmatrix x\\\xi\endbmatrix=
\bmatrix B^{-1}&-B^{-1} C\\
AB^{-1}&B-AB^{-1}C\endbmatrix\, \bmatrix y\\
\eta\endbmatrix\\\\\\\\+\,\bmatrix
B^{-1}x_0\\
AB^{-1}x_0+\eta_0\endbmatrix.
$$
When $A=0$ the phase $\Phi$
satisfies \eqref{fase} and,
consequently, the
corresponding operators are
bounded on all $M^{p,q}$.
This can also be verified by
means of the factorization
\eqref{factorization} (with
$A=0$). Indeed the continuity
of the operators
$M_{\eta_0}$, $T_{x_0}$ and
$D_B$ is easily seen, whereas
that of the Fourier
multiplier $\Fur^{-1}
U_C\Fur$ was shown, e.g., in
\cite[Lemma
2.1]{grochenig-heil}.\par
 On
the other hand, generally the
metaplectic operators are not bounded
on $M^{p,q}$ if $p\not=q$. An example
is given by the Fourier transform
itself (see \cite{feichtinger80}). An
example which instead falls in the
class of FIOs considered here is the
following one.
\begin{proposition}\label{contro}
The multiplication
$U_{I_{d}}$ is unbounded on
$M^{p,q}$, for every $1\leq
p,q\leq\infty$, with
$p\not=q$.
\end{proposition}
\begin{proof} We have $U_{I_d}f(x)=e^{\pi i
|x|^2}f(x)$. For $\lambda>0$,
we consider the
one-parameter family of
functions $f(x)=e^{-\pi
\lambda |x|^2}\in \mpq$, so
that
$\hat{f}(\o)=\lambda^{-d/2}e^{-\pi
(1/\lambda)  |\o|^2}$. For
every $1\leq p,q\leq\infty$,
by \cite[Lemma 5.3]{CN3}, we
have
$$\|f\|_{\mpq}=\|\hat{f}\|_{W(\cF L^p,L^q)}\asymp \frac{(\lambda+1)
^{d(\frac1p-\frac12)}}{\lambda^{\frac
{d}{2q}}(\lambda^2+\lambda)^{\frac
d2(\frac1p-\frac1q)}}.
$$
Since $Uf(x) e^{-\pi
(\lambda-i)|x|^2}$, so that
$\widehat{U
f}(\eta)=(\lambda-i)^{-d/2}e^{-\pi
(1/(\lambda-i)|\o|^2}$, the
same formula as above yields
$$\|Uf\|_{\mpq}\asymp \frac{[(\lambda+1)^2+1]^{\frac d2(\frac1p-\frac12)}
}{\lambda^{\frac
{d}{2q}}(\lambda^2+\lambda+1)^{\frac
d2(\frac1p-\frac1q)}}.
$$
As $\lambda\rightarrow 0$, we
have
$$\|Uf\|_{\mpq}\asymp\lambda^{-\frac {d}{2q}},\quad
\|f\|_{\mpq}\asymp\lambda^{-\frac
{d}{2q}-\frac
d2(\frac1p-\frac1q)}$$ so
that, if we assume
 $\|Uf\|_{\mpq}\leq
C\|f\|_{\mpq}$, then
$1/p-1/q\geq 0$, that is
$p\leq q$.\par Moreover, the
same argument applies to the
adjoint operator $U^\ast f
(x)= e^{-\pi i|x|^2}
f(x)$.\par Now we show that
$p=q$. By contradiction, if
$U$ were bounded on $\mpq$,
with $p<q$, its adjoint $U^*$
would satisfy
\[
\|U^\ast f\|_{M^{p',q'}}\leq
C\|f\|_{M^{p',q'}},\quad
\forall
f\in\mathcal{S}(\R^d),
\]
with $q'<p'$, which is a
contradiction to what just
proved.
\end{proof}

 \end{document}